\documentclass[onefignum,onetabnum]{siamart251216}
\usepackage{amsmath}
\usepackage{amssymb}%
\usepackage{subcaption}
\usepackage[T1]{fontenc}
\usepackage{lmodern}
\usepackage{cleveref}
\usepackage{tikz}
\usetikzlibrary{calc}
\usetikzlibrary{positioning}
\usepackage{csquotes}
\usepackage{textgreek}
\usepackage{pgfplots}
\usepackage{pgfplotstable}
\usepackage{siunitx}
\pgfplotsset{compat=1.18}

\begin{document}
\title{\textbf{A Theory of Relaxation-Based Algebraic Multigrid}\thanks{This version dated \today.\funding{This research is supported by European Union's Horizon Marie Skłodowska-Curie Actions Doctoral Networks programme, under Grant Agreement No.\ 101072344, project {AQTIVATE} (Advanced computing, QuanTum algorIthms and data-driVen Approaches for Science, Technology and Engineering) and by the Deutsche Forschungsgemeinschaft (DFG, German Research Foundation)–--Project-ID 531152215–--CRC 1701.}}}

\author{%
   Rayan Moussa\thanks{Faculty of Mathematics and Natural Sciences, Bergische Universiat\"at Wuppertal, Wuppertal 42097, Germany. \newline (moussa@uni-wuppertal.de, kkahl@uni-wuppertal.de)}
   \and
   Karsten Kahl\footnotemark[2]
  }

\date{\today}

\maketitle

\begin{abstract}
Algebraic multigrid (AMG) methods derive their optimal efficiency from the interplay between a relaxation process and a corresponding coarse grid correction. In many standard formulations, relaxation and coarse-graining are analyzed and treated as largely separate of one another. Here we propose an alternative theoretical approach centered entirely on the relaxation process, which exposes its fundamental role in the coarse-graining of the fine-scale problem. By treating the relaxation of the error as a dynamical system and applying a dimensional-reduction procedure analogous to the Mori-Zwanzig-Nakajima formalism, we derive exact expressions for the coarse-level equations and the interpolation operations, as well as a natural way of computing complementary transfer operators. We illustrate the unifying nature of this framework by recovering several well-known results for general non-symmetric systems, including ideal and optimal restriction and interpolation, as well as the limiting case of exact elimination. We further emphasize the pivotal importance of compatible-relaxation and identify dynamical corrections that naturally arise in our theory, which have the potential to enhance the convergence, robustness, and adaptivity of future algebraic multigrid methods.

\end{abstract}

\begin{keywords}
{algebraic multigrid, two-grid theory, nonsymmetric, Mori-Zwanzig, compatible-relaxation, non-Markovian memory effects, ideal/optimal restriction and interpolation, energy minimization}
\end{keywords}

\begin{AMS}
65F10, 65N22, 65N55
\end{AMS}


\section{Introduction}
Multigrid methods were originally conceived \cite{fedorenko1964speed} as a highly efficient way to overcome the fundamental flaw of classical iterative solvers: as problem size grows and conditioning worsens, classical relaxation methods such as Jacobi or Gauss-Seidel exhibit drastically deteriorating convergence rates \cite{GOODMAN199561,MGMC89,Wolff:1989wq}. Multigrid tackles this by employing a sequence of increasingly smaller problems, each a coarser representation of the one before, which allows for the elimination of the slow-to-converge parts of the computations in a far less expensive manner once the iterative method sufficiently slows down.

This multi-level philosophy was initially applied in purely geometrical settings, where coarser representations were geometrically-motivated rediscretizations of the underlying linear systems. Indeed, it has been shown that an optimal $\mathcal{O}(n)$ complexity can be achieved by applying this approach on a variety of such problems \cite{Bran1977,Bran1986, hensonmccormick2000, HackTrot1982, KahlKint2018, TrotOostSchu2001}.

The framework of \textit{Algebraic Multigrid} (AMG) \cite{RugeStue1987} was introduced to systematically construct coarse representations without any explicit dependence on geometry or discretization mode, relying entirely on the algebraic properties of the system matrix. Employing ideas such as strength-of-connection \cite{BranBranKahlLivs2015b, OlsoSchrTumi2010, RugeStue1987, Stue1983, XuZika2017}, smoothed aggregation \cite{BrezVass2011, FromKahlKrieLedeRott2014, MandBrezVane1999, Nota2010,  RottFromKahlKrieLede2012}, energy-minimization \cite{MandBrezVane1999, OlsoSchrTumi2011, WanChanSmit1999, XuZika2004}, spectral AMGe \cite{CharFalgHensJoneMantMcCoRugeVass2003, EmdeVass2001, JoneVass2001} and Adaptive or Bootstrap AMG \cite{BranBranKahlLivs2011, BranCao2022, BranKahl2014, BrezKeteMantMcCoParkRuge2012, KahlRott2018, MantMcCoParkRuge2010}, AMG methods have been extremely successful at efficiently solving high-dimensional symmetric positive-definite (SPD) systems from elliptic and parabolic PDEs. For more general non-SPD and indefinite problems, AMG faces numerous challenges which require novel generalizations to the framework of conventional AMG. In particular, the currently standard approach would be to employ the Petrov-Galerkin formulation, with techniques such as approximate ideal restriction (AIR) and reduction-based AMG \cite{AliBranKahlKrzySchrSout2024, FromKahlMacLZikaBran2010}.

Given a linear system
\begin{equation}\label{eq:linearsystem}
    {A}x={b},
\end{equation}  
with $x,b\in \mathbb{C}^n$ and $A\in \mathbb{C}^{n\times n}$ non-singular, the algebraic multigrid process can be summarized as follows. Starting from an initial approximation $x^{(0)}\in\mathbb{C}^n$, one applies $k=\mathcal{O}(1)$ iterations of a linear stationary relaxation method
\begin{equation}\label{eq:xrelaxation}
x^{(\ell+1)} = x^{(\ell)} + M r^{(\ell)},
\end{equation}
where $r^{(\ell)} = b - A x^{(\ell)}$ is the residual at iteration $\ell$, and the preconditioning matrix $M \in \mathbb{C}^{n \times n}$ determines the relaxation scheme. After $k$ steps, it is assumed that the relaxation has effectively reduced those error components well represented on the fine level, leaving a remaining error $e^{(k)}$ that can no longer be efficiently damped by relaxation alone.

The remaining error is decomposed by splitting the $n$ degrees of freedom (DOFs) into $n_c$ coarse variables and $n_f = n - n_c$ fine variables, with corresponding components $e^{(k)}_\sigma \in \mathbb{C}^{n_c}$ and $e^{(k)}_\phi \in \mathbb{C}^{n_f}$. Restriction and prolongation operators $R \in \mathbb{C}^{n_c \times n}$ and $P \in \mathbb{C}^{n \times n_c}$ are introduced to define the Petrov-Galerkin coarse operator
\begin{equation}\label{eq:RAP}
A_\sigma = R A P 
\end{equation}
and coarse residual $r^{(k)}_\sigma = R r^{(k)}$. One then solves the reduced system
\begin{equation}\label{eq:petrovgalerkin}
A_\sigma \, \epsilon^{(k)}_\sigma = r^{(k)}_\sigma
\end{equation}
for an approximation $\epsilon^{(k)}_\sigma \in \mathbb{C}^{n_c}$ of the coarse error, prolongates back to obtain $\epsilon^{(k)} = P \epsilon^{(k)}_\sigma$, and updates the iterate via $x^{(0)} \longleftarrow x^{(k)} + \epsilon^{(k)}$. The procedure is repeated until a convergence criterion is met.

The central task of any standard AMG method is to construct appropriate transfer operators $R$ and $P$. Classical AMG uses strength-of-connection heuristics to directly examine matrix entries of $A$. Smoothed Aggregation forms aggregates and relaxes them to capture the near-kernel of $A$. Energy-minimization techniques solve local optimization problems, while spectral AMGe builds operators from generalized local eigenvalue problems \cite{CharFalgHensJoneMantMcCoRugeVass2003, EmdeVass2001, JoneVass2001}. Adaptive and bootstrap AMG approaches iteratively improve the operators using test-vectors extracted from computed residuals \cite{BranBranKahlLivs2011, BranCao2022, BranKahl2014, BrezKeteMantMcCoParkRuge2012, KahlRott2018, MantMcCoParkRuge2010}. 

At the heart of virtually all standard formulations lies the following idea. Since the error $e^{(k)}$ satisfies the residual equation
    \begin{equation}\label{eq:residualequation}
        Ae^{(k)} = r^{(k)},
    \end{equation}
one can apply a coarse-graining procedure on this equation to produce a reduced problem for the coarse components $e^{(k)}_{\sigma}$. Although this allows construction of a vast landscape of successful AMG methods, it is a specific choice that does not explicitly relate the coarse-grained system to the relaxation method used. The error $e^{(k)}$ is generated entirely from $e^{(0)} = x - x^{(0)}$ via the relaxation equation\footnote{This equation derives from~\cref{eq:xrelaxation} by substituting $x^{(\ell)}$ in terms of $e^{(\ell)}$.}
    \begin{equation}\label{eq:erelaxation}
        e^{(\ell + 1)} = T e^{(\ell)}, \quad 0 \leq \ell \leq k-1,
    \end{equation}
with the \textit{error propagator} $T := I-\widehat{A}$ and the \textit{relaxation operator} $\widehat{A} := MA$.
    
    To illustrate the issue with the conventional coarse-graining strategy, consider fixing the initial iterate $x^{(0)}$ and the propagator $T$ while allowing $A$ to vary\footnote{The operator $M$ varies in a way that keeps $MA$ invariant.}. The error relaxation equation and its evolution are then unchanged, so the sequence $\{e^{(\ell)}\}$ is identical for all $\ell \in \mathcal{T}:=\{1,\dots,k\}$. Consequently, the quantities we aim to approximate, i.e., the coarse and fine components of the error and the relationship between them, remain exactly the same, and any truly relaxation-based formalism should be invariant under this change. In conventional coarse-graining, however, both the operator to be coarsened and the restricted residual depend explicitly on $A$, so the coarse operator, restricted residual, and transfer operators all vary even though the error itself does not. This reveals a redundancy in standard approaches: AMG methods built solely around the algebraic properties of $A$ alter their multigrid hierarchy even when the underlying multilevel problem is unchanged.

    Although conventional AMG has functioned effectively despite this redundancy, owing mainly to the fact that $\widehat{A}$ and $A$ are often not too different\footnote{Since $M$ is typically well-conditioned (e.g., Richardson or Jacobi), $A$ and $\widehat{A}$ largely share a common near-kernel \cite{nelson2025characterizationnearnullerrorcomponents}.}, the presence of this fundamental redundancy motivates rethinking the foundations upon which coarse-grid operators and interpolation schemes are built. A more principled approach focused on the actual dynamics of the relaxation process is needed, particularly because multigrid employs coarse-graining to approximate the evolving structure of the error that relaxation fails to eliminate, not merely to approximate static properties of the original system.

    In this paper, we address the construction of a general multigrid formalism derived entirely from the relaxation equation, which enjoys the invariance property outlined above. This is achieved by coupling relaxation and coarse-graining so that all components of the multigrid hierarchy depend exclusively on quantities intrinsic to relaxation: the error propagator $T$ and the relaxation history $\{x^{(\ell)} \mid \ell \in \mathcal{T}\}$.

    In \cref{sec:Coarse-graining}, we introduce the notation and foundational concepts underlying an invariant formulation of the coarse-graining process. This framework reveals the implicit dependence of the coarse DOFs on unresolved fine-scale interactions, which appear as \emph{dynamical memory effects} directly related to the notion of \emph{compatible-relaxation}. We present explicit expressions for the modified coarse equations and interpolation relations, which constitute the fundamental equations of the proposed framework.

    In \cref{sec:IdealizedTwoLevelConstructions}, we demonstrate how the relaxation-based equations can be used to recover both classical and novel AMG constructions in an idealized two-level setting, classified by the type and amount of information retained, in close analogy with classifications in non-equilibrium statistical mechanics \cite{chorin2005problemreductionrenormalizationmemory}.

    Finally, in \cref{sec:TransferOperators}, we construct a natural iterative update rule for relaxation-based transfer operators of increasing complementarity with the relaxation, and connect this update flow to smoothed aggregation, energy minimization and ideal/optimal interpolation and restriction operators.

\section{Problem Reduction}\label{sec:Coarse-graining}

A fundamental principle that every multigrid method must respect is \textit{locality}: achieving $\mathcal{O}(n)$ complexity requires that each step cost no more than $\mathcal{O}(1)$ computations per DOF, so any given variable can interact with only very few others. This has important implications for the relaxation, interpolation, and coarse equations.

\begin{figure}
    \begin{center}
        \resizebox{\linewidth}{!}{
            \includegraphics{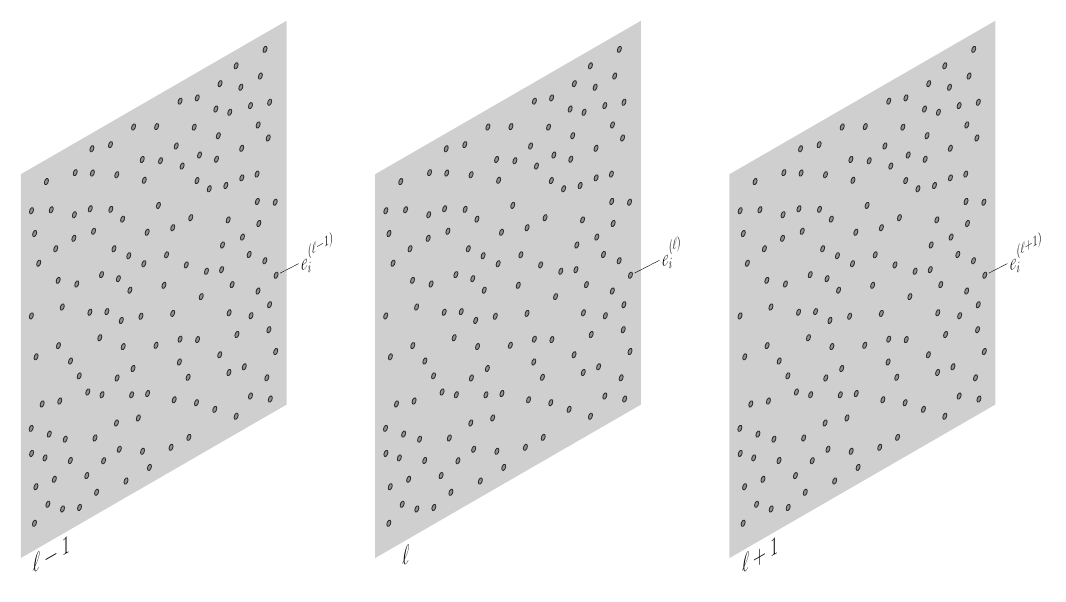}
        }
    \end{center}
        \caption{Sample layers at time-slices $\ell-1$, $\ell$, and $\ell+1$ of the spatio-temporal structure $\mathcal{N}\times\mathcal{T}$, on which the degrees of freedom $e_i^{(\ell)}$ of the discrete-time dynamical process are defined and interconnected.}
        \label{fig:spatio-temporal-graph}
\end{figure}

First and foremost, relaxation must be a local process. The most natural interpretation of equation~\eqref{eq:erelaxation} is that of a discrete-time dynamical system \cite{GOODMAN199561, MGMC89}. Describing the error vectors $e^{(0)}, \ldots, e^{(k)}$ in a fixed basis $\{g^{i} \in \mathbb{C}^{n} \mid i\in \mathcal{N}\}$, where $\mathcal{N}=\{1,\dots, n\}$, as
\begin{equation}\label{eq:errordecompositioncomponents}
    e^{(\ell)} = \sum_{i\, =\, 1}^{n}e_{i}^{(\ell)}g^{i},
\end{equation}
the DOFs of this discrete-time dynamical system can be understood as $n$ spatial points living on $k$ temporal layers, a structure denoted by $\mathcal{N}\times \mathcal{T}$ and illustrated in~\cref{fig:spatio-temporal-graph}.\footnote{This spatio-temporal structure should not be confused with discretizations of a physical spacetime domain.}
Furthermore, introducing the dual basis $\{ g_i \in \mathbb{C}^n \mid i \in \mathcal{N}\}$ satisfying the orthonormality and completeness relations
\begin{equation}\label{duality}
    g_i^\top g^{j} = \delta_i^{j},
    \qquad
    \sum_{i\,=\,1}^{n} g^{i} g_i^\top = I,
\end{equation}
and using $e_i^{(\ell)} = g_i^\top e^{(\ell)}$, the iteration equation~\cref{eq:erelaxation} becomes
\begin{equation*}
    e_i^{(\ell+1)} = \sum_{j\,=\,1}^{n} T_i^{j} \, e_j^{(\ell)},
    \quad
    T_i^{j} = g_i^\top T g^{j}.
\end{equation*}
Under the locality assumption, $T$ is sparse ($T_i^j \neq 0$ only for few indices $j$ given any $i$), and the relaxation equation takes the local component-wise form
\begin{equation}\label{releqcomponents}
    e_i^{(\ell+1)} = \sum_{j\,\in\, \mathcal{N}_{i}} T_i^{j} \, e_j^{(\ell)},
\end{equation}
where $\mathcal{N}_{i} = \{j\in \mathcal{N} \mid T_i^{j}\neq 0\}$ is the set of nearest neighbors of node $i$. This leads to the directed graph
\[
G_{T}
:
=
\Bigl(
\{(i,\ell)\in \mathcal{N}\times \mathcal{T}\},
\;\{((i,\ell),(j,\ell+1)) \mid T_{i}^{j}\neq 0\}
\Bigr),
\]
which encodes the forward-in-time propagation of information during relaxation. 

\begin{figure}
    \begin{center}
        \resizebox{\linewidth}{!}{
            \includegraphics{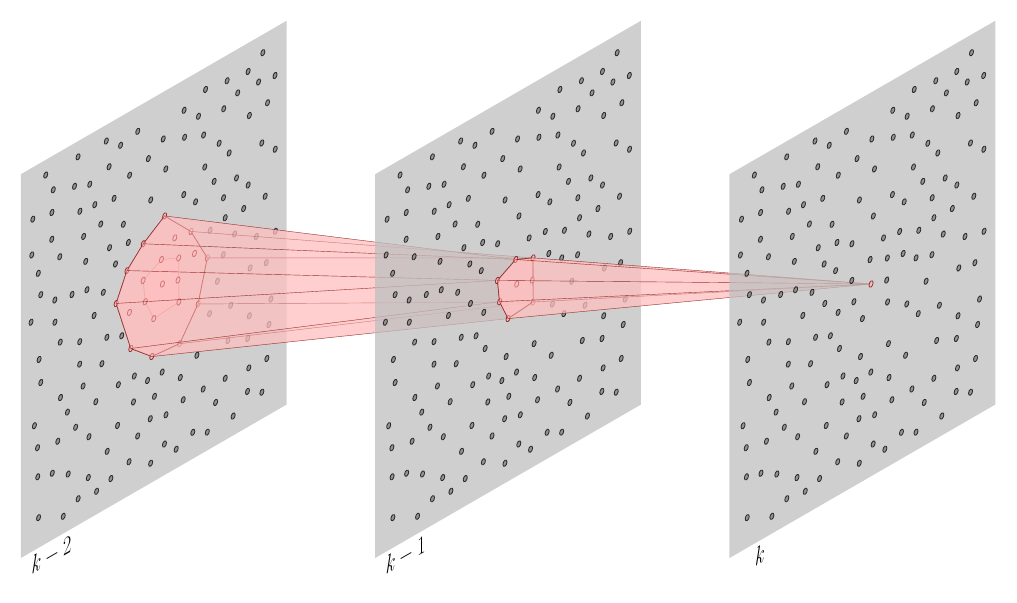}
        }%
    \end{center}
        \caption{Illustration of causal neighborhoods of a $(i,k)$, showing the nodes of $G_{T}$ which it depends on. Note, that variables on the same time slice are causally unrelated by the relaxation process since information cannot propagate within the same time-slice.} 
        \label{fig:influence-zones-backward}
\end{figure}

\begin{figure}
    \begin{center}
        \resizebox{\linewidth}{!}{
            \includegraphics{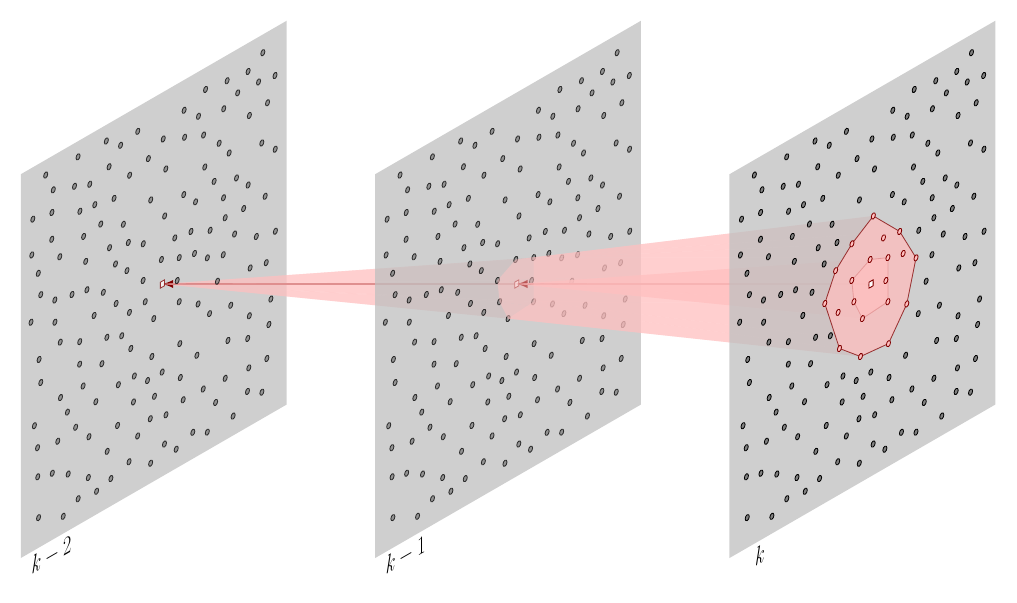}
        }%
    \end{center}
        \caption{Given $e_{i}^{(k)}$ on time slice $k$ the shift relation~\cref{eq:shift} yields $e_{i}^{(\ell)}$ for $\ell < k$, depicted by the arrows pointing backwards in time. Application of a compatible relaxation then propagates this information back to time $k$. Effectively resulting in an interpolation relation between $e_{i}^{(k)}$ and $e_{j}^{(k)}$ for all $j$ in the cone of influence of $i$.}
        \label{fig:influence-zones-interpolation}
\end{figure}

The graph's forward direction reflects that each node receives information only from earlier time-slices, so nodes on the same time-slice are causally unrelated. However, the error at a node $(j,k)$ is related to the error at $(j,\ell)$ (with $\ell \leq k$) via
\begin{equation}\label{eq:shift}
    e_{j}^{(\ell)} = e_{j}^{(k)} + x_{j}^{(k)} - x_{j}^{(\ell)},
\end{equation}
where the differences $x_{j}^{(k)} - x_{j}^{(\ell)}$ are known from the relaxation history.\footnote{This follows immediately from $e^{(\ell)} = x - x^{(\ell)}$.} Thus, once the error is known at coarse nodes on the latest time slice $k$, it is known on coarse nodes on all time slices. In combination with a forward re-propagation of information, in a manner compatible with the original relaxation process, this can reveal an effective relation between fine and coarse nodes on the same slice as illustrated in~\cref{fig:influence-zones-interpolation}.

To be more precise we denote the sets of coarse and fine DOFs by
\begin{equation*}
    \mathcal{C} \subset \mathcal{N} \text{\ and\ }  \mathcal{F} = \mathcal{N} \setminus \mathcal{C}
\end{equation*}
with $|\mathcal{C}|=n_{c}$ and $|\mathcal{F}|=n_{f}$, and set $\mathcal{C}_{i} = \mathcal{C} \cap \mathcal{N}_{i}$, $\mathcal{F}_{i} = \mathcal{F} \cap \mathcal{N}_{i}$. Splitting the sum in~\cref{releqcomponents} at the latest iteration into coarse and fine neighbor contributions,
\begin{equation*}
    e^{(k)}_{i} = \sum_{j_{1}\, \in\, \mathcal{N}_{i}} {T}_{i}^{j_{1}} e_{j_{1}}^{(k-1)} = \sum_{j_{1}\, \in\, \mathcal{C}_{i}} {T}_{i}^{j_{1}} e_{j_{0}}^{(k-1)} + \sum_{j_{1}\, \in\, \mathcal{F}_{i}} {T}_{i}^{j_{1}} e_{j_{0}}^{(k-1)}.
\end{equation*}
The first sum can then be evaluated using~\cref{eq:shift} since error values at coarse nodes are assumed known. The second sum accounts for information propagated through fine nodes on the preceding slice, which can be further decomposed into components originating from coarse nodes at $k-2$ and from fine nodes at $k-2$:
\begin{equation}\label{eq:interpolation1}
    e^{(k)}_{i} = \sum_{j_{1} \,\in\, \mathcal{C}_{i}} {T}_{i}^{j_{1}} e_{j_{1}}^{(k-1)} + \sum_{j_{1} \,\in\, \mathcal{F}_{i}}\sum_{j_{2} \,\in\, \mathcal{C}_{j_{1}}} {T}_{i}^{j_{1}} {T}_{j_{1}}^{j_{2}}e_{j_{2}}^{(k-2)} + \sum_{j_{1} \,\in\, \mathcal{F}_{i}}\sum_{j_{2} \,\in\, \mathcal{F}_{j_{1}}} {T}_{i}^{j_{1}} {T}_{j_{1}}^{j_{2}}e_{j_{2}}^{(k-2)}.
\end{equation}
Repeating this decomposition recursively, one systematically extracts all contributions that originate at a coarse node on time slice $\ell$ and pass exclusively through fine points until reaching DOF $i$ on time slice $k$. Thus we obtain 
\begin{align}\label{eq:interpolation2}
    e^{(k)}_{i}
    ={}& 
    \sum_{j_{0} \,\in\, \mathcal{C}_{i}} T_{i}^{j_{1}} \, e_{j_{1}}^{(k-1)}
    \;+\; \cdots 
    \;+\;
    \sum_{j_{1} \,\in\, \mathcal{F}_{i}}
    \cdots
    \sum_{j_{k} \,\in\, \mathcal{C}_{j_{k-1}}}
    T_{i}^{j_{1}} \cdots T_{j_{k-1}}^{j_{k}}
    \, e_{j_{k}}^{(0)}
    \nonumber \\[0.3em]
    &\;+\;
    \sum_{j_{1} \,\in\, \mathcal{F}_{i}}
    \cdots
    \sum_{j_{k} \,\in\, \mathcal{F}_{j_{k-1}}}
    T_{i}^{j_{1}} \cdots T_{j_{k-1}}^{j_{k}}
    \, e_{j_{k}}^{(0)} .
\end{align}

The final term represents information that originates at fine nodes on time slice $0$ and propagates along paths of length $k$ that pass exclusively through fine nodes at every intermediate step. This term cannot be accounted for by the knowledge of the error on $\mathcal{C}$.  In order to for the coarse-level error approximation to accurately represent the fine-node error via interpolation, this term needs to be negligible, i.e., information propagating exclusively along fine-only paths must decay sufficiently fast. This requirement is precisely the principle of \emph{compatible-relaxation}, a notion extensively studied and used in standard algebraic multigrid literature (see, e.g., \cite{Bran2000,BranCaoKahlFalgHu2018,BranFalg2010,BranZika2007,DAmbVass2013,Livn2004}).\footnote{The notion of compatible relaxation used here corresponds to the general habituated variant found in the literature.}
Formally, the effective coupling
\begin{equation*}
\widetilde{T}^{j_k}_{i}
:= \sum_{j_{1} \,\in\, \mathcal{F}_{i}} \cdots 
  \sum_{j_{k-1} \,\in\, \mathcal{F}_{j_k}}
  T_{i}^{j_{1}} \cdots T_{j_{k-1}}^{j_k}
\end{equation*}
must be sufficiently small; achieving this decay rests solely on an appropriate choice of basis vectors and duals as well as their split into $\mathcal{C}$ and $\mathcal{F}$.
\begin{figure}
    \begin{center}
        \resizebox{\linewidth}{!}{
            \includegraphics[scale=1]{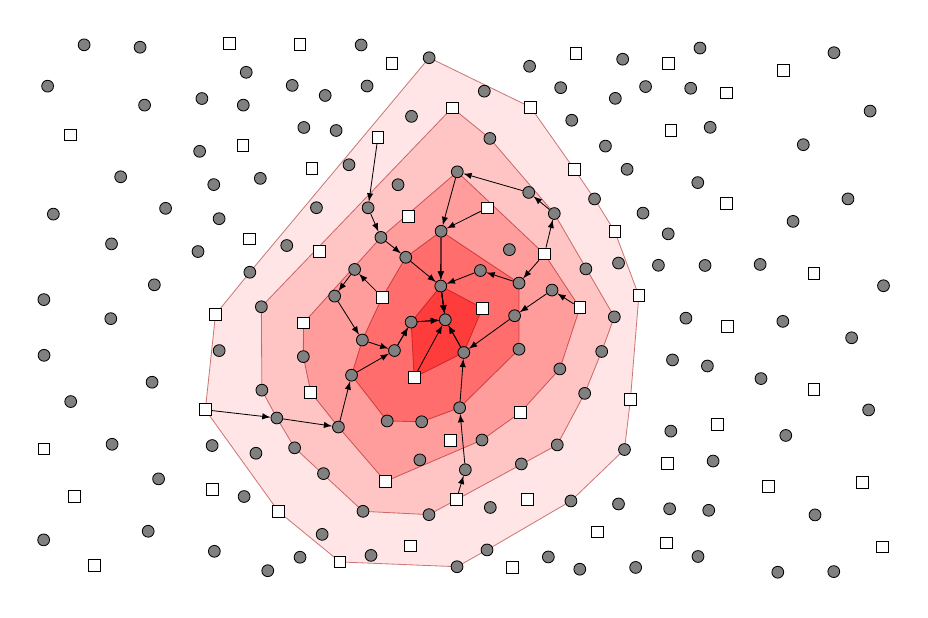}
        }%
    \end{center}
        \caption{Illustration of some of the terms appearing in~\cref{eq:interpolation}. Each term is a path in $G_{T}$ originating at coarse nodes (squares) and ending at the fine node of interest while passing only through fine points (circles) on its way. The colors of the level sets represent the size of the contribution provided by the corresponding terms associated with coarse DOFs as suggested by the principle of compatible relaxation.}
        \label{fig:influence-zones-paths}
\end{figure}

Under the assumption that compatible-relaxation is well satisfied, the final term in~\eqref{eq:interpolation2} can be neglected. The resulting approximation $\epsilon_{i}^{(k)}$ of $e_{i}^{(k)}$ then fulfills the interpolation relation
\begin{equation}\label{eq:interpolation}
\epsilon_{i}^{(k)} = \sum_{\ell\, =\, 0}^{k-1}\sum_{j\,\in\, \mathcal{C}^{(\ell+1)}_{(i,k)}}  \sum_{j_{\ell}\,\in\, \mathcal{F}_{j}}\widetilde{T}_{i}^{j_{\ell}} T_{j_{\ell}}^{j} \epsilon_{j}^{(k-\ell-1)},
\end{equation}
where $\mathcal{C}^{(\ell+1)}_{(i,k)}$ denotes the set of spatial indices corresponding to coarse nodes of $G_{T}$ which are connected to $(i,k)$ via fine-paths of length $\ell + 1$.
Introducing the effective interpolation weights
\begin{equation*}
W_{ij}^{(\ell)} := \sum_{j_{\ell}\, \in\, \mathcal{F}_{j}} \widetilde{T}_{i}^{j_{\ell}} T_{j_{\ell}}^{j},
\end{equation*}
the interpolation equation assumes the form
\begin{equation}\label{eq:winterpolation}
\epsilon_{i}^{(k)} = \sum_{\ell\, =\, 0}^{k-1}\sum_{j\in \mathcal{C}^{(\ell+1)}_{(i,k)}} W_{ij}^{(\ell)} \epsilon_{j}^{(k-\ell-1)},
\end{equation}
highlighting the temporal non-locality induced by relaxation: the fine error at iteration $k$ depends on coarse error values from earlier iterations, with weights determined by effective long-range fine-fine couplings as illustrated in~\cref{fig:influence-zones-paths}.
\subsection{Matrix Form}

Having exposed interpolation in component form and explained its dynamical origin we now introduce a simplified matrix representation which allows us to expose its connection to standard formulations of AMG. 

Given a partition of the DOFs into a coarse set $\mathcal{C}$ and a fine set $\mathcal{F}$, we collect the corresponding error components into vectors $e_{\sigma}^{(\ell)}\in\mathbb{C}^{n_c}$ and $e_{\phi}^{(\ell)}\in\mathbb{C}^{n_f}$, and assemble the associated basis vectors as columns of the coarse and fine prolongation operators
\begin{equation}\label{ProlongationDef}
    P = \begin{bmatrix}
        p^{1} \mid \cdots \mid p^{n_{c}}
    \end{bmatrix}
    \in \mathbb{C}^{n \times n_{c}},
    \qquad
    Q = \begin{bmatrix}
        q^{1} \mid \cdots \mid q^{n_{f}}
    \end{bmatrix}
    \in \mathbb{C}^{n \times n_{f}},
\end{equation}
and their duals as rows of
\begin{equation}\label{ReductionDef}
    P^{\dagger} = 
    \begin{bmatrix}
        p_{1} \mid \cdots \mid p_{n_{c}}
    \end{bmatrix}^\top
    \in \mathbb{C}^{n_{c} \times n},
    \qquad
    Q^{\dagger} = 
    \begin{bmatrix}
        q_{1} \mid \cdots \mid q_{n_{f}}
    \end{bmatrix}^\top
    \in \mathbb{C}^{n_{f} \times n}.
\end{equation}
The component-wise relations $e_{\sigma}^{(\ell)} = P^{\dagger} e^{(\ell)}$ and $e_{\phi}^{(\ell)} = Q^{\dagger} e^{(\ell)}$ hold immediately. The duality relation~\eqref{duality} yields the normalization conditions
\begin{equation*}
    P^{\dagger} P = I,
    \qquad
    Q^{\dagger} Q = I,
\end{equation*}
the orthogonality relations $Q^{\dagger} P = 0$, $P^{\dagger} Q = 0$, and the completeness relation
\begin{equation*}
    PP^{\dagger} + QQ^{\dagger} = I,
\end{equation*}
so that $PP^{\dagger}$ and $QQ^{\dagger}$ are oblique projections onto the coarse and fine subspaces, and the error admits the decomposition
\begin{equation}\label{eq:errordecomp}
    e^{(\ell)} = P e^{(\ell)}_{\sigma} + Q e^{(\ell)}_{\phi}.
\end{equation}
The transition elements $T_{i}^{j}$ take the matrix form
\begin{equation*}
  T_{i}^{ j}=\begin{cases}
    \left(Q^{\dagger}TQ \right)_{ij}, & i,j \in \mathcal{F},\\
    \left(Q^{\dagger}TP \right)_{ij}, & i \in \mathcal{F},\; j\in \mathcal{C}.
  \end{cases}
\end{equation*}
The effective transition elements from coarse node $j$ to fine node $i$ over $\ell$ fine-only steps satisfy
\begin{equation*}
    \sum_{j_{1}\, \in\, \mathcal{F}_{i}}\cdots\sum_{j\, \in\, \mathcal{C}_{j_{\ell-1}}} {T}_{i}^{j_{1}}\cdots{T}_{j_{\ell-1}}^{j} = \left((Q^{\dagger}TQ)^{\ell}Q^{\dagger}TP \right)_{ij},
\end{equation*}
and the generalized interpolation law~\cref{eq:interpolation2} can be written in matrix form as
\begin{equation}\label{eq:interpolationmatform}
        e_{\phi}^{(k)}= \sum_{\ell\, =\, 0}^{k-1}(Q^{\dagger}TQ)^{\ell}Q^{\dagger}TP \, e_{\sigma}^{(k-\ell-1)} + (Q^{\dagger}TQ)^{k}e^{(0)}_{\phi}.
\end{equation}
This decomposition explicitly exposes the information content of the fine-grid error: contributions propagated from the coarse grid over successive steps of compatible relaxation through the repeated application of $Q^{\dagger}TQ$, plus a residual term confined to the fine grid that depends only on $e_\phi^{(0)}$.

In this representation the special role of the compatible relaxation propagator $Q^{\dagger}TQ$ is revealed as it not only governs the propagation of information from coarse DOFs to fine degrees of freedom, but also the the final term, which is the only contribution that cannot be obtained from coarse-grid information as $e^{(0)}_\phi$ is arbitrary and unknown.

\begin{figure}
    \centering
    \includegraphics[width=1\linewidth]{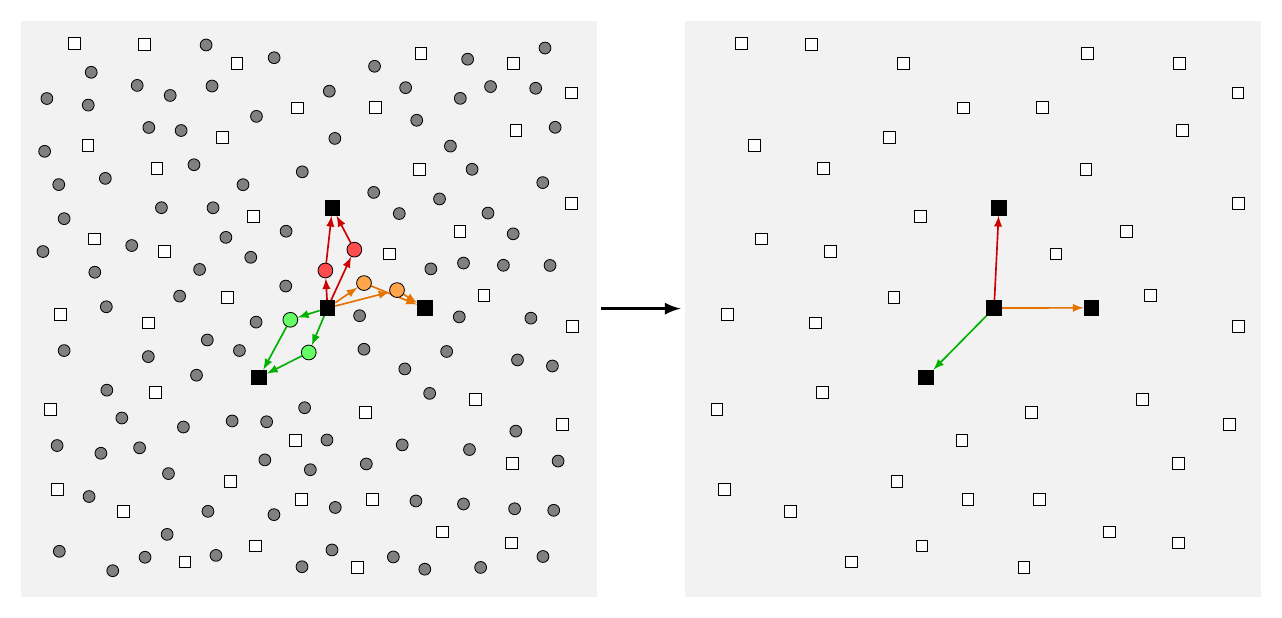}
    \caption{Transition from a fine-level description of relaxation to a coarse-level description. The coarse nodes (squares) interact with each other through intermediate fine nodes (circles). The effective interaction between the coarse nodes in the coarse representation is then equivalent to the original interaction between the coarse nodes in the original fine representation through fine-only paths which originally connect them. Here only a subset of these paths are drawn.}
    \label{fig:coarse-graining}
\end{figure}

Assuming compatible-relaxation holds, i.e., $(Q^{\dagger}TQ)^k \rightarrow 0$ rapidly, and given an approximation $\epsilon^{(k)}_{\sigma}$ to the coarse-node error, the last term of~\cref{eq:interpolationmatform} may be dropped, giving the interpolation law
\begin{equation}\label{eq:winterpolationmatform}
    \epsilon_{\phi}^{(k)} = \sum_{\ell\, =\, 0}^{k-1}W^{(\ell+1)}\epsilon^{(k-\ell-1)}_{\sigma},
\end{equation}
where the {generalized interpolation weights} matrix is
\begin{equation*}
    W^{(\ell+1)} := (Q^{\dagger} T Q)^{\ell} Q^{\dagger} \widehat{A} P, \quad \ell \geq 0.
\end{equation*}
Using~\cref{eq:shift} at the coarse points, the coarse-grid error at earlier iterations is approximated by
\begin{equation}\label{eq:shiftmatform}
    \epsilon_{\sigma}^{(\ell)}
    =
    \epsilon_{\sigma}^{(k)} + x_{\sigma}^{(k)} - x_{\sigma}^{(\ell)}.
\end{equation}
Then combining the decomposition~\cref{eq:errordecomp} with the interpolation equation~\cref{eq:winterpolationmatform} gives
\begin{equation}\label{eq:eapprox}
    \epsilon^{(k)}
    =
    P^{(0)} \epsilon_{\sigma}^{(k)}
    +
    P^{(1)} \epsilon_{\sigma}^{(k-1)}
    + \cdots +
    P^{(k)} \epsilon_{\sigma}^{(0)},
 \end{equation}
with the {generalized prolongation operators}
\begin{equation}\label{eq:Pell}
    P^{(0)} := P,
    \qquad
    P^{(\ell)} := Q W^{(\ell)}, \quad \ell \geq 1.
\end{equation}
The coarse-grained relaxation equation can then be written as
\begin{equation}\label{eq:morizwanzigmatform}
    \epsilon^{(k+1)}_{\sigma}
    =
    \sum_{\ell\, =\, 0}^{k-1}
    T^{(\ell)} \, \epsilon_{\sigma}^{(k-\ell)},
\end{equation}
where the effective coarse-level propagators
\begin{equation*}
    T^{(\ell)} := P^{\dagger} T P^{(\ell)}, \qquad \ell \geq 0,
\end{equation*}
encode memory effects introduced by fine-grid variables whose elimination results in a non-local coarse evolution, as depicted in~\cref{fig:coarse-graining}. This structure closely parallels the reduced dynamics obtained in projection-based coarse-graining frameworks such as the Mori-Zwanzig-Nakajima formalism \cite{chorin2005problemreductionrenormalizationmemory,etde_6362584, LinTianLiveAngh2021,Mori:1965oqj,Nakajima:1958pnl, Zwanzig:1960gvu}, and clarifies the role of fine-grid relaxation as a source of memory in the coarse-grid evolution.

\subsection{The Coarse Equation}
To obtain a coarse equation solvable for the coarse DOFs, we rewrite the relaxation equation~\cref{eq:erelaxation} as $e^{(k+1)} - e^{(k)} = -\widehat{A} e^{(k)}$. Since $e^{(k+1)} - e^{(k)} = x^{(k)} - x^{(k+1)}$, this gives
\begin{equation*}
    \widehat{A} e^{(k)} = \widehat{r}^{(k)},
    \qquad
    \widehat{r}^{(k)} := x^{(k+1)} - x^{(k)},
\end{equation*}
a preconditioned form of the standard residual equation~\cref{eq:residualequation}. Applying a restriction operator $R\in \mathbb{C}^{n_{c}\times n}$ and substituting the generalized interpolation relation~\cref{eq:interpolationmatform} yields the \emph{exact} coarse residual equation
\begin{equation}\label{eq:exactcoarseequation}
    \sum_{\ell\, =\, 0}^{k}
    A_{\sigma}^{(\ell)} \, e_{\sigma}^{(k-\ell)}
    =
    \widehat{r}_{\sigma}^{(k)} + \eta^{(k)},
\end{equation}
where $\widehat{r}_{\sigma}^{(k)} := R \widehat{r}^{(k)}$, the sequence of coarse operators is
\begin{equation}\label{eq:RAPell}
    A_{\sigma}^{(\ell)} := R \widehat{A} P^{(\ell)},
\end{equation}
and the noise term\footnote{This terminology is motivated by the Mori-Zwanzig formalism, where an analogous term represents the interaction of a coarse-grained system with its eliminated degrees of freedom.}
\begin{equation}\label{eq:noise}
    \eta^{(k)} := - R \widehat{A} Q (Q^{\dagger} T Q)^{k} e_{\phi}^{(0)}
\end{equation}
collects contributions from unresolved fine-grid error components. If $\eta^{(k)}$ is made small and thus negligible, we obtain the generalized coarse residual equation
\begin{equation}\label{eq:genreseq}
    \sum_{\ell\, =\, 0}^{k}
    A_{\sigma}^{(\ell)} \, \epsilon_{\sigma}^{(k-\ell)}
    =
    \widehat{r}_{\sigma}^{(k)},
\end{equation}
which, together with~\cref{eq:shiftmatform}, completely determines $\epsilon^{(k)}_{\sigma}$. Like~\cref{eq:morizwanzigmatform}, this relation is non-Markovian, i.e., it depends not only on the latest error, but rather on the full history of the error evolution.

For the choice $R = P^{\dagger}$, for which $\widehat{r}_{\sigma}^{(k)} = \epsilon_{\sigma}^{(k)} - \epsilon_{\sigma}^{(k+1)}$, one can verify directly that~\cref{eq:genreseq} is equivalent to the coarse-grained relaxation equation~\cref{eq:morizwanzigmatform}, confirming that it is a valid coarse-level description usable to solve for the required coarse error variables.

\section{Idealized Two-Level Constructions}\label{sec:IdealizedTwoLevelConstructions}
\subsection{Markovian}\label{sec:MarkovianConstruction} 
The coarse equations~\eqref{eq:genreseq} and interpolation relations~\cref{eq:interpolation} are inherently non-Markovian. To recover a Markovian formulation, we drop any memory terms, which results in the interpolation relation
\[
\epsilon^{(k)} = P \epsilon_{c}^{(k)}
\] and coarse equation given by the Petrov-Galerkin formulation
\begin{equation}\label{eq:petrovgalerkinequation}
    R \widehat{A} P \, \epsilon_{\sigma}^{(k)} = \widehat{r}_{\sigma}^{(k)}.
\end{equation} Equation~\eqref{eq:petrovgalerkinequation} has the same structure as the standard Petrov-Galerkin equation~\cref{eq:petrovgalerkin}, with $A$ replaced by $\widehat{A}$ and the residual by $\widehat{r}_{\sigma}^{(k)}$. Defining $\widehat{R} := R M^{-1}$ (assuming $M$ invertible) reduces~\eqref{eq:RAQ} to the standard Petrov-Galerkin form with coarse operator $A_{\sigma} = \widehat{R} A P$ and residual $\widehat{R} r^{(k)}$.

Solving the coarse equation for $e_{\sigma}^{(k)}$ reconstructs the coarse error components exactly, leaving
\begin{equation}\label{eq:Qe}
  e^{(k)} - \epsilon^{(k)} = Q Q^{\dagger} T^{k} e^{(0)}.
\end{equation}
as the remaining error. The coarse-level solve can also be represented by the projection $QQ^{\dagger} = I - P (R \widehat{A} P)^{-1} R \widehat{A}$, giving the two-grid error propagator
\begin{equation}\label{eq:markovianpropagator}
    E_{\mathrm{TG}}^{(k)}
    =
    \left(
        I - P (R \widehat{A} P)^{-1} R \widehat{A}
    \right)
    \left(
        I - \widehat{A}
    \right)^{k},
\end{equation}
which takes a well-known form which is frequently encountered in the context of standard two-level AMG theory;\footnote{Albeit being expressed in terms of $\widehat{A}$ rather than $A$, as required by the invariance principle.} see, e.g., \cite{FalgVassZika2005, XuZhan2018, 392c305c-406c-3726-9762-23620c1a010a, AliBranKahlKrzySchrSout2025}.

Taking a close look at the terms that are dropped in the Markovian formulation reveals specific conditions in which the memory terms disappear naturally. On one hand, we observe that due to
\[
P^{(\ell)} = \begin{cases}
    P & \text{for\ } \ell = 0\\
    \left(Q^{\dagger}TQ\right)^{\ell-1}Q^{\dagger}\widehat{A}P & \text{for\ } \ell >0 
\end{cases}
\] the memory terms $\ell > 0$ vanish if and only if $Q^{\dagger}\widehat{A}P = 0.$ This corresponds to transfer operators which are $\widehat{A}$ (or $T$) invariant.\footnote{This is fulfilled for optimal interpolation as discussed in~\cref{sec:optimalinterpolation}.} With regards to the coarse equation on the other hand we obtain from definitions~\cref{eq:RAPell} and~\cref{eq:Pell}, that $A_{\sigma}^{(\ell)} = 0$ for all $\ell \geq 1$ if
\begin{equation}\label{eq:RAQ}
    R \widehat{A} Q = 0.
\end{equation}
That is the subspaces spanned by the rows of $R$ and the columns of $Q$ are $\widehat{A}$-orthogonal. Under this condition, the coarse equation reduces to~\cref{eq:petrovgalerkinequation}, the noise term vanishes identically ($\eta^{(k)} = 0$), and the solution is exact: $\epsilon_{\sigma}^{(k)} = e_{\sigma}^{(k)}$. 
\subsection{Semi-Markovian}
The Markovian formulation discards all terms in equations~\cref{eq:genreseq} and~\cref{eq:interpolationmatform} which depend on the history of error iterates. Suppose we have $R\widehat{A}Q = 0$ so that we end up with a Markovian coarse equation~\cref{eq:petrovgalerkinequation} to solve for $e_{\sigma}^{(k)}$, but then use the shift relation~\cref{eq:shiftmatform} to reconstruct the full history $\{e_{\sigma}^{(\ell)} \mid 1 \leq \ell \leq k\}$:
\begin{equation}\label{eq:shiftmatformmarkoviancoarsesolve}
    e_{\sigma}^{(\ell)}
    =
    (R\widehat{A}P)^{-1}\widehat{r}^{(k)}_{\sigma}
    + x_{\sigma}^{(k)} - x_{\sigma}^{(\ell)}.
\end{equation}
Substituting into the non-Markovian interpolation equation~\cref{eq:eapprox} gives the affine correction
\begin{equation}\label{eq:winterpolationmatformaffine}
    \epsilon^{(k)}
    =
    \sum_{\ell\, =\, 0}^{k-1}
    P^{(\ell)}
    \left(
        (R\widehat{A}P)^{-1}\widehat{r}^{(k)}_{\sigma}
        + x_{\sigma}^{(k)} - x_{\sigma}^{(k-\ell)}
    \right).
\end{equation}
Since the Markovian coarse solve is exact and all terms in the generalized interpolation are accounted for, the only remaining contribution is the unresolved compatible-relaxation component. Therefore, the two-grid error propagator, as can be easily shown, reduces to
\begin{equation*}
    E_{\text{TG}}^{(k)}
    =
    Q (Q^{\dagger} T Q)^{k} Q^{\dagger},
\end{equation*}
showing that combining an exact Petrov-Galerkin coarse correction with memory-inclusive interpolation yields a two-level method whose convergence is entirely governed by habituated compatible-relaxation.

\subsection{Non-Markovian}\label{sec:nonmarkovian}
Relaxing the orthogonality constraint~\cref{eq:RAQ}, the non-Markovian coarse equation~\cref{eq:genreseq} must be rewritten in a form solvable for $e_{\sigma}^{(k)}$. Using the shift relation~\cref{eq:shift} to eliminate dependence on previous iterates, it becomes
\begin{equation}\label{eq:nonmarkoviancoarseequation}
    A_{\sigma}\,\epsilon_{\sigma}^{(k)}
    =
    \widehat{r}_{\sigma}^{(k)} + s_{\sigma}^{(k)},
\end{equation}
where the generalized coarse operator and memory correction are
\begin{equation*}
    A_{\sigma}
    :=
    \sum_{\ell\, =\, 0}^{k} A_{\sigma}^{(\ell)},
    \qquad
    s_{\sigma}^{(k)}
    :=
    -
    \sum_{\ell\, =\, 0}^{k} A_{\sigma}^{(k-\ell)}
    \left( x_{\sigma}^{(k)} - x_{\sigma}^{(\ell)} \right).
\end{equation*}
The solution of~\cref{eq:nonmarkoviancoarseequation} differs from the exact solution of~\cref{eq:exactcoarseequation} by $A_{\sigma}^{-1} \eta^{(k)}$, and the full remaining error is
\begin{equation*}
    e^{(k)} - \epsilon^{(k)}
    =
    \left(
        I - P' A_{\sigma}^{-1} R \widehat{A}
    \right)
    Q (Q^{\dagger} T Q)^{k} Q^{\dagger} e^{(0)},
\end{equation*}
where the effective prolongation operator is
\begin{equation}\label{eq:P'}
    P' := \sum_{\ell\, =\, 0}^{k} P^{(\ell)}.
\end{equation}
The two-level error propagator takes the form
\begin{equation}\label{eq:nonmarkovianpropagator}
    E_{\mathrm{TG}}^{(k)}
    =
    \left[
        I - P' (R \widehat{A} P')^{-1} R \widehat{A}
    \right]
    Q \bigl[ Q^{\dagger} (I - \widehat{A}) Q \bigr]^{k} Q^{\dagger},
\end{equation}
which resembles a Petrov-Galerkin construction with $P'$ in place of $P$, where the compatible-relaxation propagator takes the role of the error relaxation propagator $T$. Imposing~\cref{eq:RAQ} recovers the semi-Markovian propagator.

An important aspect of both the semi-Markovian and non-Markovian schemes is that, unlike the Markovian approach~\cref{eq:markovianpropagator}, they can be rendered exact by an appropriate choice of $Q$ and $Q^{\dagger}$. Specifically, if the dual vectors are chosen $T$-orthogonal to the basis vectors, i.e., $q_i^\top T q^{j} = 0$ for all $i,j$, then $Q^{\dagger} T Q = 0$ and the non-Markovian propagator~\eqref{eq:nonmarkovianpropagator} vanishes, yielding an exact two-level method.

\subsection{Exact}

If we wish to obtain an exact solution for an arbitrary coarse-fine split without the special basis required in the non-Markovian case, we must explicitly resolve the compatible-relaxation term using only coarse values, thus necessarily violating the locality constraint. Using
\begin{equation*}
    e_{\phi}^{(0)} = e_{\phi}^{(k)} + x_{\phi}^{(k)} - x_{\phi}^{(0)}
\end{equation*}
and substituting into~\cref{eq:winterpolationmatform} gives an equation that can be inverted to yield the exact interpolation law
\begin{equation}\label{eq:twolevelexactinterpol}
    e_{\phi}^{(k)}
    =
    \sum_{\ell\, =\, 0}^{k-1} \tilde{W}^{(\ell+1)} e_{\sigma}^{(k-\ell-1)}
    +
    \xi_{\phi}^{(k)},
\end{equation}
with modified weights and affine correction
\begin{align*}
    \tilde{W}^{(\ell)}
    & =
    \left[I - (Q^{\dagger}TQ)^{k}\right]^{-1} W^{(\ell)} \text{\ and}
    \\\nonumber
    \xi_{\phi}^{(k)}
    & =
    \left[I - (Q^{\dagger}TQ)^{k}\right]^{-1}
    (Q^{\dagger}TQ)^{k}
    \left(x_{\phi}^{(k)} - x_{\phi}^{(0)}\right).
\end{align*}
Using this exact reconstruction, decomposing and restricting the relaxation equation yields the exact coarse equation
\begin{equation}\label{eq:twolevelexactcoarse}
    \widetilde{A}_{\sigma} e_{\sigma}^{(k)}
    =
    \widehat{r}_{\sigma}^{(k)} + \widetilde{s}_{\sigma}^{(k)},
\end{equation}
with the coarse operator $\widetilde{A}_{\sigma}$ and the correction to the residual $\widetilde{s}_{\sigma}^{(k)}$ defined as
\begin{equation*}
    \widetilde{A}_{\sigma} = \sum_{\ell\,=\, 0}^{k} \widetilde{A}_{\sigma}^{(\ell)} \text{\ and\ }
    \widetilde{s}_{\sigma}^{(k)}
    =
    -\sum_{\ell\,=\,0}^{k} \widetilde{A}_{\sigma}^{(k-\ell)}
      \left(x_{\sigma}^{(k)} - x_{\sigma}^{(\ell)}\right)
    - R\widehat{A}Q\,\xi_{\phi}^{(k)}.
\end{equation*}
The corresponding iteration-dependent transfer operators are then given by $P^{(0)} = P$ and $R^{(0)} = R$. For $\ell > 1$ we find
\begin{equation*}
P^{(\ell)} :=
Q \widetilde{W}^{(\ell)} \text{\ and\ }
R^{(\ell)} := - R\widehat{A}Q
\left[I-(Q^{\dagger}TQ)^{k}\right]^{-1}
(Q^{\dagger}TQ)^{\ell-1}Q^{\dagger}.
\end{equation*}
Using the geometric series identity \[[I-(Q^{\dagger}TQ)^{k}]^{-1}\sum_{\ell\,=\,0}^{k-1}(Q^{\dagger}TQ)^{\ell} = (Q^{\dagger}\widehat{A}Q)^{-1},\] the generalized coarse operator simplifies to the Petrov-Galerkin form
\begin{equation*}
    A_{\sigma} = R\widehat{A}\widetilde{P} = \widetilde{R}\widehat{A}P
\end{equation*}
in terms of effective transfer operators
\begin{equation*}
    \widetilde{P}  
    =
    \left(I - Q(Q^{\dagger}\widehat{A}Q)^{-1}Q^{\dagger}\widehat{A}\right)P \text{\ and\ } 
    \widetilde{R}
    =
    R\!\left(I - \widehat{A}Q(Q^{\dagger}\widehat{A}Q)^{-1}Q^{\dagger}\right)
\end{equation*}
or equivalently as $A_{\sigma} = R\widehat{\mathcal{A}}P$ with $\widehat{\mathcal{A}} = \widehat{A} - \widehat{A}Q(Q^{\dagger}\widehat{A}Q)^{-1}Q^{\dagger}\widehat{A}$.

Solving~\cref{eq:twolevelexactcoarse} and interpolating via~\cref{eq:twolevelexactinterpol} yields an exact two-level method valid for \emph{any} choice of transfer operators. Although the global inversions involved render this approach impractical as an algorithm, it is conceptually important as it traces a continuous hierarchy from simple relaxation, through progressively richer coarse-grained constructions, to an exact direct solver. Conventional multigrid is the lowest-order truncation; higher-order corrections incorporate additional memory effects and recover more fine-scale information. In the limit where all this information is preserved, the AMG method converges in one iteration for any problem and any transfer operators, provided at least one relaxation step has been taken.

\section{Transfer Operators}\label{sec:TransferOperators}
The non-Markovian constructions developed above suggest a systematic viewpoint: rather than prescribing transfer operators a priori, one can \emph{derive} them directly from the relaxation dynamics. We now show how the derived formulas can be used to iteratively construct prolongation matrices\footnote{It can also be shown that an analogous iterative procedure can be derived for the remaining transfer operators.} of increasing quality, and connect this construction to classical AMG notions.

\subsection{Ideal}
The non-Markovian coarse operator $A_{\sigma}$ is equivalent to a Petrov-Galerkin coarse operator for the effective prolongation
\begin{equation}\label{eq:Pprime}
    P' = P + Q\sum_{\ell\, =\, 0}^{k-1}(Q^{\dagger}TQ)^{\ell}Q^{\dagger}T P.
\end{equation}
As the Petrov-Galerkin coarse operator built from $P'$ produces more accurate coarse equations than that based solely on $P$ (fewer terms are dropped from~\cref{eq:exactcoarseequation}), it is natural to expect that $P'$ represents the subspace which relaxation struggles to eliminate more accurately than the original $P$.

Starting from the simplest possible choice of transfer operators containing only canonical unit basis vectors, sorted so that
\begin{equation}\label{eq:P0Q0}
    P_{0} = \begin{bmatrix}
        I \\ 0
    \end{bmatrix}, \quad
    Q_{0} = \begin{bmatrix}
        0 \\ I
    \end{bmatrix}, \quad
    P_{0}^{\dagger} = \begin{bmatrix}
        I & 0
    \end{bmatrix}, \quad
    Q_{0}^{\dagger} = \begin{bmatrix}
        0 & I
    \end{bmatrix},
\end{equation}
and arranging $T$ in block form $T = \bigl[\begin{smallmatrix} T_{\mathit{cc}} & T_{\mathit{cf}} \\ T_{\mathit{fc}} & T_{\mathit{ff}} \end{smallmatrix}\bigr]$, the improved prolongation operator $P'$ evaluates to
\begin{equation}\label{eq:P1}
    P_{1}=\begin{bmatrix}
        I \\ W
    \end{bmatrix} \text{\ with\ }
    W = \sum_{\ell\, =\, 0}^{k-1}T_{\mathit{ff}}^{\ell}T_{\mathit{fc}} = -\sum_{\ell\, =\, 0}^{k-1}(I-\widehat{A}_{\mathit{ff}})^{\ell}\widehat{A}_{\mathit{fc}}.
\end{equation}
In the limit of a large number of relaxation steps, and given that the spectral radius of $T_{\mathit{ff}}$ is less than one, i.e., relaxatio is convergent, this matrix converges to a generalization of the \emph{ideal interpolation} matrix
\begin{equation*}
    W_{\text{ideal}} = \lim_{k\rightarrow\infty}W = -\widehat{A}_{\mathit{ff}}^{-1}\widehat{A}_{\mathit{fc}}.
\end{equation*}
In case the preconditioner $M$ is block upper-triangular so that $\widehat{A}_{\mathit{fc}}=M_{\mathit{ff}}A_{\mathit{fc}}$ and $\widehat{A}_{\mathit{ff}}=M_{\mathit{ff}}A_{\mathit{ff}}$, then (assuming $M_{\mathit{ff}}$ invertible) this reduces to the standard form $-A_{\mathit{ff}}^{-1}A_{\mathit{fc}}$, expressed in terms of the original operator $A$ as in conventional AMG.

The complementary dual operator satisfying $Q_1^\dagger P_1 = 0$ is
\begin{equation}\label{eq:Qd1}
    Q^{\dagger}_{1} = \begin{bmatrix}
        -W & I
    \end{bmatrix},
\end{equation}
and the ideal restriction, derived from condition~\cref{eq:RAQ}, is (up to left multiplication by an arbitrary $n_c\times n_c$ matrix)
\begin{equation*}
    R =
    \begin{bmatrix}
        I & -\widehat{A}_{\mathit{cf}}\widehat{A}_{\mathit{ff}}^{-1}
    \end{bmatrix}.
\end{equation*}
Rewriting this expression in terms of $\widehat{R} = RM^{-1}$ recovers the standard ideal restriction $\widehat{R} = [I\ -A_{\mathit{cf}}A_{\mathit{ff}}^{-1}]$ in terms of $A$.

Both ideal operators involve $\widehat{A}_{\mathit{ff}}^{-1}$ (or ${A}_{\mathit{ff}}^{-1}$ in the standard case), which is generally dense and computationally impractical since it violates the sparsity constraint. In practice, this inversion can be replaced by localized or sparse approximations, yielding \emph{approximate ideal interpolation} and \emph{approximate ideal restriction} as the basis of many practical AMG schemes \cite{AliBranKahlKrzySchrSout2024,XuZhan2018}.

\subsection{The Flow Equation}
Expression~\cref{eq:Pprime} yields an update rule for a recursive improvement of prolongation matrices:
\begin{equation}\label{eq:flow}
    P_{\tau+1}
    = P_{\tau}
    + Q_{\tau}\sum_{\ell\, =\, 0}^{k-1}(Q_{\tau}^{\dagger}TQ_{\tau})^{\ell}Q_{\tau}^{\dagger}T P_{\tau},
\end{equation}
where $\tau \in \mathbb{N}$ is the flow time and the iteration is initialized with a set of initial tentative transfer operators.\footnote{Of course, they do not have to assume the simple forms suggested in the set of equations (\ref{eq:P0Q0})} Introducing the projection $\mathcal{Q}_{\tau} := Q_{\tau}Q_{\tau}^{\dagger}$, the update simplifies to
\begin{equation*}
    P_{\tau+1}
    = \sum_{\ell\, =\, 0}^{k}(\mathcal{Q}_{\tau}T)^{\ell} P_{\tau},
\end{equation*}
which can be read as $k$ steps of relaxation applied to $P_\tau$ using the propagator $\mathcal{Q}_{\tau}T$ rather than $T$ itself, the latter being used instead in methods such as Smoothed Aggregation \cite{BrezFalgMacLMantMcCoRuge2004, SterMantMcCoMillPearRugeSand2010, BrezVass2011}. 

The flow is intrinsically nonlinear since the projection $\mathcal{Q}_\tau = I - P_\tau P_\tau^\dagger$ depends on the current state of the prolongation and its dual at every iteration of the flow. Furthermore, a key structural property is that the update $(P_{\tau+1} - P_\tau)$ is always $P_\tau^\dagger$-orthogonal:
$P_{\tau}^{\dagger}(P_{\tau+1}-P_{\tau}) = 0$, since $P_{\tau}^{\dagger}\mathcal{Q}_{\tau}=0$.
Consequently, $P_\tau^\dagger P_{\tau+1} = I_c$, so the dual operator can be chosen as stationary: $P_{\tau}^{\dagger} = P_0^{\dagger}$ for all $\tau$, simplifying the projection to $\mathcal{Q}_\tau = I - P_\tau P_0^\dagger$.

At a fixed point of the flow, i.e., $P_\infty = \lim_{\tau\to\infty} P_\tau$, the flow equation implies
\begin{equation}\label{eq:invariantsubspace}
    Q_{\infty}^{\dagger}T P_{\infty} = 0.
\end{equation}
Together with $Q_\infty^\dagger P_\infty = 0$, this shows that $P_\infty$ encodes a right $T$-invariant subspace, $\mathrm{col}(TP_\infty) \subseteq \mathrm{col}(P_\infty)$, while $Q_\infty^\dagger$ spans the orthogonal left $T$-invariant subspace.

In the limit $k\to\infty$,\footnote{Well approximated for finite $k$ when compatible-relaxation is fast to converge.} the flow equation becomes
\begin{equation*}
   P_{\tau+1}
   = P_{\tau}
   - Q_{\tau}(Q_{\tau}^{\dagger}\widehat{A}Q_{\tau})^{-1}
     Q_{\tau}^{\dagger}\widehat{A} P_{\tau},
\end{equation*}
which, using the projection $ \Pi_{\widehat{A}} =:  I - Q_{\tau}(Q_{\tau}^{\dagger}\widehat{A}Q_{\tau})^{-1}
     Q_{\tau}^{\dagger}\widehat{A}$, can be rewritten 
     as 
$$P_{\tau+1}
   =\Pi_{\widehat{A}} P_{\tau},$$
implying that the flow
projects $P_\tau$ onto the $\widehat{A}$-orthogonal complement of $\mathrm{col}(Q_\tau)$ at every instance of the flow.

When $\widehat{A}$ is Hermitian and $Q_\tau$ has properly orthonormalized columns, meaning $Q_{\tau}^{\dagger} = Q_{\tau}^{H}$, the projection $\Pi_{\widehat{A}}$ is $\widehat{A}$-orthogonal and the \textit{energy} of the individual columns $p^{i}_{\tau}$, defined as
\begin{equation*}
    \|p^i_\tau\|_{\widehat{A}} = p^{i\; H}_{\tau} \widehat{A} p^{i}_{\tau},
\end{equation*}
is non-increasing under the flow, i.e., 
\begin{equation*}
    \|p^i_{\tau+1}\|_{\widehat{A}} \leq \|p^i_{\tau}\|_{\widehat{A}}.
\end{equation*}
The flow thus defines an energy-decreasing evolution of the prolongation operator as a whole, fitting naturally within the framework of \emph{energy-minimization techniques} \cite{MandBrezVane1999, OlsoSchrTumi2011, WanChanSmit1999, XuZika2004}.

\subsection{Optimal}\label{sec:optimalinterpolation}
As the fixed points of the flow represent invariant subspaces of the relaxation operator $T$, and ones of minimal-energy if a corresponding energy-minimization criteria are met, the flow should converge to the so-called \emph{optimal interpolation} subspace. The eigenmodes that relax most slowly are those whose eigenvalues have the largest magnitude. Ordering the spectrum of $T$ by decreasing magnitude, we separate coarse and fine blocks:
\begin{equation*}
    \Lambda = \mathrm{diag}(\underbrace{\lambda_1,\dots,\lambda_{n_c}}_{\Lambda_c},\underbrace{\lambda_{n_c+1},\dots,\lambda_n}_{\Lambda_f}),
\end{equation*}
and the right and left eigenvector matrices $V_{R}$ and $V_{L}$, which satisfy $TV_R = V_R\Lambda$ and $V_LT = \Lambda V_L$, as $V_R = [V_R^{n_c}\ V_R^{n_f}]$ and $V_L = [(V_L^{n_c})^\top\;(V_L^{n_f})^\top]^\top$. Normalizing so that $V_L V_R = I$ yields the block biorthogonality relations:
\begin{equation}\label{eq:eigenorthog}
    V_{L}^{n_{c}}V_{R}^{n_{c}} = I_{c}, \quad
    V_{L}^{n_{f}}V_{R}^{n_{c}} = 0, \quad
    V_{L}^{n_{c}}V_{R}^{n_{f}} = 0, \quad
    V_{L}^{n_{f}}V_{R}^{n_{f}} = I_{f}.
\end{equation}

If the flow indeed converges to the invariant subspace of minimal energy, the fixed-point prolongation must span the optimal coarse eigenspace:
\begin{equation*}
    P_\infty = V_R^{n_c}\mathcal{Z}_c, \qquad P_\infty^\dagger = \mathcal{Z}_c^{-1}V_L^{n_c},
\end{equation*}
for some invertible $\mathcal{Z}_c$ determined by the initial condition on the transfer operators. The complementary operators spanning the high-energy invariant subspace are thus
\begin{equation*}
    Q_\infty = V_R^{n_f}\mathcal{Z}_f^{-1}, \qquad Q_\infty^\dagger = \mathcal{Z}_f V_L^{n_f},
\end{equation*}
and the ideal restriction operator of equation~\cref{eq:RAQ} is thus $R_\infty = \mathcal{Z}_c' V_L^{n_c}$.

With these choices, the non-Markovian formulas reduce to the semi-Markovian ones via~\cref{eq:RAQ}, and further to standard Petrov-Galerkin once the invariant subspace condition~\cref{eq:invariantsubspace} is imposed. The coarse correction then eliminates the slow-to-converge modes entirely, and the non-Markovian two-grid error propagator reduces to
\begin{equation*}
   E_{\mathrm{TG}} = V_R^{n_f}\Lambda_f^k V_L^{n_f},
\end{equation*}
meaning that the asymptotic convergence factor is governed solely by $\lambda_{n_c+1}$. The theory of optimal interpolation and restriction has been extensively studied in standard AMG literature~\cite{AliBranKahlKrzySchrSout2025, BranCaoKahlFalgHu2018, KrzySoutWimmAliBranKahl25} and the
fact that the present non-Markovian framework reduces to the standard optimal construction under the flow of transfer operators further supports its optimality. No further improvement can be achieved by exploiting information from memory terms, which is due to $Q_{\infty}^{\dagger}\widehat{A}P_{\infty} = 0$ effectively reducing the framework to the standard Markovian case as already discussed in~\cref{sec:MarkovianConstruction}.

\section{Conclusions}
We have explored the construction of the fundamental AMG equations in a two-level setting, though the framework naturally generalizes to multiple levels by performing a dynamical coarse-graining of the relaxation process rather than the more traditional coarse graining of the residual equation. This perspective is motivated by the concept of dynamical invariance of the AMG hierarchy: the hierarchy should remain fixed under perturbations of the operator provided the error relaxation dynamics are preserved.

Using a formalism closely related to the projective Mori-Zwanzig-Nakajima approach, we derived the corresponding coarse equations by systematically eliminating fine-scale variables while retaining their dynamical influence through memory terms. The resulting equations are inherently non-Markovian, reflecting the dependence of coarse variable evolution on the full relaxation history.

Rather than coarse-graining the residual equation associated with $A$, we coarse-grain the relaxation equation through the propagator $T$ and the relaxation operator $\widehat{A}$. This shift in perspective leads naturally to a hierarchy of generalized multigrid schemes, interpolating between classical AMG and increasingly more accurate higher-order corrections. As more dynamical information is retained, the method approaches an exact direct solve with correspondingly improving convergence properties.

We provided a dynamical interpretation of all terms arising in this construction, relating them to non-Markovian memory effects, compatible-relaxation, energy minimization, and ideal or optimal interpolation and restriction, demonstrating that several classical AMG ideas emerge as natural consequences of the underlying coarse-grained dynamics.

Although some exact formulations are impractical due to global inversions, they serve as a useful theoretical baseline suggesting a continuous spectrum of methods from inexpensive local approximations to increasingly more accurate but costlier corrections. Understanding this trade-off may guide the design of more robust, adaptive, and efficient multigrid algorithms. We view this framework as a foundation for developing principled strategies that balance accuracy and complexity in practical AMG solvers.

While this paper focuses on the theoretical aspects of the relaxation-based framework for AMG, a second manuscript illustrating the findings with extensive numerical tests is currently being worked on. Finding practical and efficient ways to exploit memory effects in AMG both for symmetric and non-symmetric problems is a topic of future research.

\bibliography{lib}

@article{AliBranKahlKrzySchrSout2025,
	author = {Ali, Ahsan and Brannick, James and Kahl, Karsten and Krzysik, Oliver A. and Schroder, Jacob B. and Southworth, Ben S.},
	title = {Generalized Optimal {AMG} Convergence Theory for Nonsymmetric and Indefinite Problems},
	journal = {SIAM Journal on Scientific Computing},
	volume = {0},
	number = {0},
	pages = {S89--S111},
	year = {2025},
	doi = {10.1137/24M1679288}
}

@article{AliBranKahlKrzySchrSout2024,
	author = {Ahsan Ali and James J. Brannick and Karsten Kahl and Oliver A. Krzysik and Jacob B. Schroder and Ben S. Southworth},
	title = {Constrained Local Approximate Ideal Restriction for Advection-Diffusion Problems},
	journal = {SIAM Journal on Scientific Computing},
	volume = {0},
	number = {0},
	pages = {S96–S122},
	year = {2024},
	doi = {10.1137/23M1583442}
}

@article{Bran1977,
	author = {Achi Brandt},
	title = {Multi-Level Adaptive Solutions to Boundary-Value Problems},
	journal = {Mathematics of Computation},
	volume = {31},
	number = {138},
	pages = {333–390},
	year = {1977},
	doi = {10.1090/S0025-5718-1977-0431719-X}
}

@article{Bran1986,
	author = {Achi Brandt},
	title = {Algebraic Multigrid Theory: The Symmetric Case},
	journal = {Applied Mathematics and Computation},
	volume = {19},
	number = {1–4},
	pages = {23–56},
	year = {1986},
	doi = {10.1016/0096-3003(86)90095-0}
}

@article{Bran2000,
	author = {Achi Brandt},
	title = {General highly accurate algebraic coarsening},
	journal = {Electronic Transactions on Numerical Analysis},
	volume = {10},
	number = {1},
	pages = {21},
	year = {2000}
}

@article{BranBranKahlLivs2011,
	author = {Achi Brandt and James Brannick and Karsten Kahl and Ira Livshits},
	title = {Bootstrap {AMG}},
	journal = {SIAM Journal on Scientific Computing},
	volume = {33},
	number = {2},
	pages = {612–632},
	year = {2011},
	doi = {10.1137/090752973}
}

@article{BranBranKahlLivs2015b,
	author = {Achi Brandt and James Brannick and Karsten Kahl and Ira Livshits},
	title = {Algebraic Distance for Anisotropic Diffusion Problems: Multilevel Results},
	journal = {Electronic Transactions on Numerical Analysis},
	volume = {44},
	pages = {472–496},
	year = {2015}
}

@article{BranCaoKahlFalgHu2018,
	author = {James Brannick and Fei Cao and Karsten Kahl and Robert D. Falgout and Xiaozhe Hu},
	title = {Optimal Interpolation and Compatible Relaxation in Classical Algebraic Multigrid},
	journal = {SIAM Journal on Scientific Computing},
	volume = {40},
	number = {3},
	pages = {A1473–A1493},
	year = {2018},
	doi = {10.1137/17M1123456}
}

@article{BranCao2022,
	author = {James Brannick and Shuhao Cao},
	title = {A Bootstrap Multigrid Eigensolver},
	journal = {SIAM Journal on Matrix Analysis and Applications},
	volume = {43},
	number = {4},
	pages = {1627–1657},
	year = {2022},
	doi = {10.1137/20M131151X}
}

@article{BranFalg2010,
	author = {James Brannick and Robert D. Falgout},
	title = {Compatible relaxation and coarsening in algebraic multigrid},
	journal = {SIAM Journal on Scientific Computing},
	volume = {32},
	number = {3},
	pages = {1393–1416},
	year = {2010},
	doi = {10.1137/090772216}
}

@article{BranKahl2014,
	author = {James Brannick and Karsten Kahl},
	title = {Bootstrap Algebraic Multigrid for the 2D {W}ilson {D}irac system},
	journal = {SIAM Journal on Scientific Computing},
	volume = {36},
	number = {3},
	pages = {B321–B347},
	year = {2014},
	doi = {10.1137/130934660}
}

@inproceedings{BranZika2007,
	author = {James Brannick and Ludmil T. Zikatanov},
	title = {Algebraic Multigrid Methods Based on Compatible Relaxation and Energy Minimization},
	booktitle = {Domain Decomposition Methods in Science and Engineering XVI},
	editor = {Olof B. Widlund and David E. Keyes},
	pages = {15–26},
	publisher = {Springer},
	address = {Berlin, Heidelberg},
	year = {2007},
	doi = {10.1007/978-3-540-34469-8_2}
}

@article{BrezFalgMacLMantMcCoRuge2004,
	author = {Marian Brezina and Robert D. Falgout and Scott P. MacLachlan and Thomas A. Manteuffel and Stephen F. McCormick and John W. Ruge},
	title = {Adaptive Smoothed Aggregation (\textalpha {SA}).},
	journal = {SIAM Journal on Scientific Computing},
	volume = {25},
	number = {6},
	pages = {1896–1920},
	year = {2004},
	doi = {10.1137/S1064827502418598}
}

@article{BrezKeteMantMcCoParkRuge2012,
	author = {Marian Brezina and Christian Ketelsen and Thomas A. Manteuffel and Stephen F. McCormick and Minho Park and John W. Ruge},
	title = {Relaxation-corrected bootstrap algebraic multigrid (r{BAMG})},
	journal = {Numerical Linear Algebra with Applications},
	volume = {19},
	number = {2},
	pages = {178–193},
	year = {2012},
	doi = {10.1002/nla.1821}
}

@inproceedings{BrezVass2011,
	author = {Marian Brezina and Panayot S Vassilevski},
	title = {Smoothed aggregation spectral element agglomeration AMG: SA--$\rho$ AMGe},
	booktitle = {International Conference on Large--Scale Scientific Computing},
	pages = {3–15},
	year = {2011}
}

@book{hensonmccormick2000,
	author = {Briggs, William and Henson, Van and McCormick, Steve},
	title = {A Multigrid Tutorial, 2nd Edition},
	year = {2000},
	isbn = {978-0-89871-462-3}
}

@article{CharFalgHensJoneMantMcCoRugeVass2003,
	author = {Timothy P. Chartier and Robert D. Falgout and Van Enden Henson and Jim E. Jones and Thomas A. Manteuffel and Stephen F. McCormick and John W. Ruge and Panayot S. Vassilevski},
	title = {Spectral {AMG}e (\textrho {AMG}e).},
	journal = {SIAM Journal on Scientific Computing},
	volume = {25},
	number = {1},
	pages = {1–26},
	year = {2003},
	doi = {10.1137/S106482750139892X}
}

@article{chorin2005problemreductionrenormalizationmemory,
	author = {Alexandre J. Chorin and Panagiotis Stinis},
	title = {Problem reduction, renormalization, and memory},
	journal = {Communications in Applied Mathematics and Computational Science},
	volume = {1},
	number = {1},
	pages = {1--27},
	year = {2006},
	doi = {10.2140/camcos.2006.1.1}
}

@article{DAmbVass2013,
	author = {Pasqua D’Ambra and Panayot S. Vassilevski},
	title = {Adaptive {AMG} with coarsening based on compatible weighted matching},
	journal = {Computing and Visualization in Science},
	volume = {16},
	number = {2},
	pages = {59–76},
	year = {2013},
	doi = {10.1007/s00791-014-0224-9}
}

@article{EmdeVass2001,
	author = {Henson Van Emden and Panayot S. Vassilevski},
	title = {Element-Free {AMG}e: General Algorithms for Computing Interpolation Weights in {AMG}},
	journal = {SIAM Journal on Scientific Computing},
	volume = {23},
	number = {2},
	pages = {629–650},
	year = {2001},
	doi = {10.1137/S1064827500372997}
}

@article{FalgVassZika2005,
	author = {Robert D. Falgout and Panayot S. Vassilevski and Ludmil T. Zikatanov},
	title = {On Two-Grid Convergence Estimates},
	journal = {Numerical Linear Algebra with Applications},
	volume = {12},
	number = {5--6},
	pages = {471–494},
	year = {2005},
	doi = {10.1002/nla.437}
}

@article{fedorenko1964speed,
	author = {Fedorenko, Radii Petrovich},
	title = {The speed of convergence of one iterative process},
	journal = {USSR Computational Mathematics and Mathematical Physics},
	volume = {4},
	number = {3},
	pages = {227--235},
	year = {1964}
}

@article{FromKahlKrieLedeRott2014,
	author = {Andreas Frommer and Karsten Kahl and Stefan Krieg and Björn Leder and Matthias Rottmann},
	title = {Adaptive Aggregation-Based Domain Decomposition Multigrid for the Lattice {W}ilson--{D}irac Operator},
	journal = {SIAM Journal on Scientific Computing},
	volume = {36},
	number = {4},
	pages = {A1581–A1608},
	year = {2014},
	doi = {10.1137/130919507}
}

@article{FromKahlMacLZikaBran2010,
	author = {Andreas Frommer and Karsten Kahl and Scott P. MacLachlan and Ludmil T. Zikatanov and James Brannick},
	title = {Adaptive Reduction-Based Multigrid for Nearly Singular and Highly Disordered Physical systems},
	journal = {Electronic Transactions on Numerical Analysis},
	volume = {37},
	pages = {276–295},
	year = {2010}
}

@article{GOODMAN199561,
	author = {Jonathan Goodman and Neal Madras},
	title = {Random-walk interpretations of classical iteration methods},
	journal = {Linear Algebra and its Applications},
	volume = {216},
	pages = {61--79},
	year = {1995},
	doi = {10.1016/0024-3795(93)00101-5}
}

@article{MGMC89,
	author = {Goodman, Jonathan and Sokal, Alan D.},
	title = {Multigrid {M}onte {C}arlo method. {C}onceptual foundations},
	journal = {Phys. Rev. D},
	volume = {40},
	pages = {2035--2071},
	year = {1989},
	doi = {10.1103/PhysRevD.40.2035}
}

@book{etde_6362584,
	author = {Hermann Grabert},
	title = {Projection Operator Techniques in Nonequilibrium Statistical Mechanics},
	publisher = {Springer},
	address = {Berlin, Heidelberg},
	series = {Springer Tracts in Modern Physics},
	volume = {95},
	year = {1982},
	doi = {10.1007/BFb0044591}
}

@book{HackTrot1982,
	author = {Wolfgang Hackbusch and Ulrich Trottenberg},
	title = {Multigrid Methods},
	publisher = {Springer},
	series = {Lecture Notes in Mathematics},
	year = {1982},
	isbn = {978-3-540-39544-7}
}

@article{JoneVass2001,
	author = {Jim E. Jones and Panayot S. Vassilevski},
	title = {{AMG}e Based on Element Agglomeration.},
	journal = {SIAM Journal on Scientific Computing},
	volume = {23},
	number = {1},
	pages = {109–133},
	year = {2001},
	doi = {10.1137/S1064827599361047}
}

@article{KahlKint2018,
	author = {Karsten Kahl and Nils Kintscher},
	title = {Geometric Multigrid for the Tight-Binding {H}amiltonian of Graphene.},
	journal = {SIAM Journal on Numerical Analysis},
	volume = {56},
	number = {1},
	pages = {499–519},
	year = {2018},
	doi = {10.1137/16M1102033}
}

@article{KahlRott2018,
	author = {Karsten Kahl and Matthias Rottmann},
	title = {Least Angle Regression Coarsening in Bootstrap Algebraic Multigrid},
	journal = {SIAM Journal on Scientific Computing},
	volume = {40},
	number = {6},
	pages = {A3928–A3954},
	year = {2018},
	doi = {10.1137/18M1168182}
}

@article{KrzySoutWimmAliBranKahl25,
	author = {Krzysik, Oliver A. and Southworth, Ben S. and Wimmer, Golo A. and Ali, Ahsan and Brannick, James and Kahl, Karsten},
	title = {Optimal Transfer Operators in Algebraic Two-Level Methods for Nonsymmetric and Indefinite Problems},
	year = {2025}
}

@article{LinTianLiveAngh2021,
	author = {Yen Ting Lin and Yifeng Tian and Daniel Livescu and Marian Anghel},
	title = {Data-Driven Learning for the {Mori}-{Zwanzig} Formalism: A Generalization of the {K}oopman Learning Framework},
	journal = {SIAM Journal on Applied Dynamical Systems},
	volume = {20},
	number = {4},
	pages = {2558–2601},
	year = {2021},
	doi = {10.1137/21M1401759}
}

@article{Livn2004,
	author = {O. E. Livne},
	title = {Coarsening by compatible relaxation.},
	journal = {Numerical Lin. Alg. with Applic.},
	volume = {11},
	number = {2--3},
	pages = {205–227},
	year = {2004},
	doi = {10.1002/nla.378}
}

@article{MandBrezVane1999,
	author = {J. Mandel and Marian Brezina and P. Vaněk},
	title = {Energy Optimization of Algebraic Multigrid Bases},
	journal = {Computing},
	volume = {62},
	number = {3},
	pages = {205–228},
	year = {1999},
	doi = {10.1007/s006070050022}
}

@article{MantMcCoParkRuge2010,
	author = {Thomas A. Manteuffel and Stephen F. McCormick and Minho Park and John W. Ruge},
	title = {Operator-based interpolation for bootstrap algebraic multigrid.},
	journal = {Numerical Linear Algebra with Applications},
	volume = {17},
	number = {2--3},
	pages = {519–537},
	year = {2010},
	doi = {10.1002/nla.711}
}

@article{Mori:1965oqj,
	author = {Mori, Hazime},
	title = {Transport, Collective Motion, and {B}rownian Motion},
	journal = {Prog. Theor. Phys.},
	volume = {33},
	number = {3},
	pages = {423--455},
	year = {1965},
	doi = {10.1143/PTP.33.423}
}

@article{Nakajima:1958pnl,
	author = {Nakajima, Sadao},
	title = {On Quantum Theory of Transport Phenomena: Steady Diffusion},
	journal = {Prog. Theor. Phys.},
	volume = {20},
	number = {6},
	pages = {948--959},
	year = {1958},
	doi = {10.1143/PTP.20.948}
}

@misc{nelson2025characterizationnearnullerrorcomponents,
	author = {Austen J. Nelson and Panayot S. Vassilevski},
	title = {Characterization of the near-null error components utilized in composite adaptive AMG solvers},
	year = {2025}
}

@article{Nota2010,
	author = {Yvan Notay},
	title = {An aggregation-based algebraic multigrid method},
	journal = {Electronic Transactions on Numerical Analysis},
	volume = {37},
	pages = {123–146},
	year = {2010}
}

@article{392c305c-406c-3726-9762-23620c1a010a,
	author = {YVAN NOTAY},
	title = {ANALYSIS OF TWO-GRID METHODS: THE NONNORMAL CASE},
	journal = {Mathematics of Computation},
	volume = {89},
	number = {322},
	pages = {807--827},
	year = {2020},
	doi = {10.1090/mcom/3478}
}

@article{OlsoSchrTumi2010,
	author = {Luke N. Olson and Jacob B. Schroder and Raymond S. Tuminaro},
	title = {A new perspective on strength measures in algebraic multigrid.},
	journal = {Numerical Linear Algebra with Applications},
	volume = {17},
	number = {4},
	pages = {713–733},
	year = {2010},
	doi = {10.1002/nla.669}
}

@article{OlsoSchrTumi2011,
	author = {Luke N. Olson and Jacob B. Schroder and Raymond S. Tuminaro},
	title = {A General Interpolation Strategy for Algebraic Multigrid Using Energy Minimization},
	journal = {SIAM Journal on Scientific Computing},
	volume = {33},
	number = {2},
	pages = {966–991},
	year = {2011},
	doi = {10.1137/100803031}
}

@inproceedings{RottFromKahlKrieLede2012,
	author = {Matthias Rottmann and Andreas Frommer and Karsten Kahl and Stefan Krieg and Björn Leder},
	title = {Aggregation-based Multilevel Methods for Lattice {QCD}},
	booktitle = {Proceedings of XXIX International Symposium on Lattice Field Theory — PoS(Lattice 2011)},
	pages = {046},
	year = {2012},
	doi = {10.22323/1.139.0046}
}

@inbook{RugeStue1987,
	author = {John W. Ruge and Klaus Stüben},
	title = {Algebraic Multigrid},
	booktitle = {Multigrid Methods},
	editor = {Stephen F. McCormick},
	chapter = {4},
	pages = {73–130},
	publisher = {Society for Industrial and Applied Mathematics},
	series = {Frontiers in Applied Mathematics},
	year = {1987},
	doi = {10.1137/1.9781611971057.ch4}
}

@article{SterMantMcCoMillPearRugeSand2010,
	author = {Hans De Sterck and Thomas A. Manteuffel and Stephen F. McCormick and Killian Miller and J. Pearson and John W. Ruge and Geoffrey Sanders},
	title = {Smoothed Aggregation Multigrid for Markov Chains.},
	journal = {SIAM Journal on Scientific Computing},
	volume = {32},
	number = {1},
	pages = {40–61},
	year = {2010},
	doi = {10.1137/080719157}
}

@article{Stue1983,
	author = {Klaus Stüben},
	title = {Algebraic Multigrid ({AMG}): Experiences and Comparisons},
	journal = {Applied Mathematics and Computation},
	volume = {13},
	number = {3--4},
	pages = {419–451},
	year = {1983},
	doi = {10.1016/0096-3003(83)90023-1}
}

@book{TrotOostSchu2001,
	author = {Ulrich Trottenberg and Cornelis W. Oosterlee and Anton Schüller},
	title = {Multigrid},
	publisher = {Academic Press},
	address = {San Diego [u.a.]},
	series = {Texts in Applied Mathematics. Bd.},
	volume = {33},
	year = {2001},
	note = {With contributions by A. Brandt, P. Oswald and K. Stüben},
	isbn = {0-12-701070-X}
}

@article{WanChanSmit1999,
	author = {W. Wan and T. Chan and B. Smith},
	title = {An Energy-minimizing Interpolation for Robust Multigrid Methods},
	journal = {SIAM Journal on Scientific Computing},
	volume = {21},
	number = {4},
	pages = {1632–1649},
	year = {1999},
	doi = {10.1137/S1064827598334277}
}

@article{Wolff:1989wq,
	author = {Wolff, Ulli},
	title = {Critical Slowing Down},
	journal = {Nucl. Phys. B Proc. Suppl.},
	volume = {17},
	pages = {93--102},
	year = {1990},
	doi = {10.1016/0920-5632(90)90224-I}
}

@article{XuZika2004,
	author = {Jinchao Xu and Ludmil T. Zikatanov},
	title = {On an energy minimizing basis for algebraic multigrid methods},
	journal = {Computing and Visualization in Science},
	volume = {7},
	number = {3},
	pages = {121–127},
	year = {2004},
	doi = {10.1007/s00791-004-0147-y}
}

@article{XuZika2017,
	author = {Jinchao Xu and Ludmil T. Zikatanov},
	title = {Algebraic multigrid methods},
	journal = {Acta Numerica},
	volume = {26},
	pages = {591–721},
	year = {2017},
	doi = {10.1017/S0962492917000083}
}

@article{XuZhan2018,
	author = {Xuefeng Xu and Chen-Song Zhang},
	title = {On the ideal interpolation operator in algebraic multigrid methods},
	journal = {SIAM Journal on Numerical Analysis},
	volume = {56},
	number = {3},
	pages = {1693–1710},
	year = {2018},
	doi = {10.1137/17M1138822}
}

@article{Zwanzig:1960gvu,
	author = {Zwanzig, Robert},
	title = {Ensemble Method in the Theory of Irreversibility},
	journal = {J. Chem. Phys.},
	volume = {33},
	number = {5},
	pages = {1338},
	year = {1960},
	doi = {10.1063/1.1731409}
}
\end{document}